\documentclass[12pt]{article}
\usepackage{amssymb,amsmath,amscd,amsthm}
\usepackage{hyperref}
\usepackage{graphicx,color,float}

\newtheorem{theorem}{Theorem}[section]

\newtheorem{lemma}[theorem]{Lemma}
\newtheorem{proposition}[theorem]{Proposition}

\theoremstyle{remark}
\newtheorem{remark}[theorem]{Remark}

\newcommand{\thmref}[1]{Theorem~\ref{#1}}
\newcommand{\propref}[1]{Proposition~\ref{#1}}

\newcommand{\lemref}[1]{Lemma~\ref{#1}}

\newcommand{\remref}[1]{Remark~\ref{#1}}

\def\eps{{\varepsilon}}
\def\cL{\mathcal{L}}
\def\cD{\mathcal{D}}
\def\cF{\mathcal{F}}
\def\bu{\mathbf{u}}
\def\R{\mathbb{R}}
\def\E{\mathbb{E}}
\def\P{\mathbb{P}}
\def\Cint{\mathcal{C}^{\rm int}_\eps}
\def\Cext{\mathcal{C}^{\rm ext}_\eps}

\begin{document}

\date{}
\title{Distribution of particles near the front in supercritical branching Brownian motion with compactly supported branching}

\author{Pratima Hebbar\footnote{Department of Mathematics,
Grinnell College, Grinnell, IA.
Email: hebbarpr@grinnell.edu}, Leonid Koralov\footnote{Department of Mathematics,
University of Maryland,
College Park, MD. Email: koralov@umd.edu}}

\maketitle

\begin{abstract}
We study a $d$-dimensional supercritical branching Brownian motion whose branching and killing rates are continuous and compactly supported. Let $\lambda_0>0$ be the principal eigenvalue of the associated Schr\"odinger operator, $\psi$ its ground state, and $W_\infty$ the limit of the additive martingale $e^{-\lambda_0 t}\sum_i \psi(X^{(i)}_t)$. We consider the number of particles in a translate of a fixed bounded set located at a distance $\sqrt{\lambda_0/2}\, t - \frac{d-1}{2\sqrt{2\lambda_0}}\ln t - b(t)$ from the origin in a fixed direction $\mathbf{u}$. Using sharp global asymptotics for parabolic equations with a compactly supported potential, we show that this number converges in distribution, together with all its moments, to a mixed Poisson random variable whose random intensity is an explicit multiple of $W_\infty$ depending on the limit of $b(t)$ and $\mathbf{u}$. Deeper inside the front (when $b(t) \rightarrow \infty$), the number of particles normalized by its expectation converges to $W_\infty/\psi(x)$. Results of this type were previously known in dimension one and, for test sets of macroscopic angular size, in dimension two; here we treat test sets of bounded size in every dimension.
\end{abstract}



\section{Introduction}

We consider a collection of particles that perform independent $d$-dimensional Brownian motions. Let $\cD$ be a compact set in $\R^d$ and let $\alpha$ and $\beta$ be continuous non-negative functions supported in $\cD$, with $\alpha\not\equiv 0$. A particle located at $x$ branches into two particles at rate $\alpha(x)$ and is annihilated at rate $\beta(x)$. Each newly created particle starts at the location of its parent and then evolves independently of all the other particles. Let $v = \alpha - \beta$.

Let $n^x_t(U)$ be the number of particles in a bounded nonempty open set $U \subset \R^d$ at time $t$, assuming that at time zero there was a single particle located at $x$. The density of particles $\rho_1(t,x,y)$ satisfies
\begin{equation} \label{dens}
\frac{\partial\rho_1(t,x,y)}{\partial t} = \frac{1}{2} \Delta_x \rho_1(t,x,y) + v(x)\rho_1(t,x,y),\qquad \rho_1(0,x,y) = \delta_y(x),
\end{equation}
where $y$ is treated as a parameter. The large-time behavior of $n^x_t(U)$ depends crucially on whether the operator
\begin{equation} \label{eqL}
\cL u  = \frac{1}{2} \Delta u  + v  u
\end{equation}
has a positive eigenvalue. We consider the supercritical case, where it does. Since $v$ is continuous and compactly supported, $\cL$ is self-adjoint on $L^2(\R^d)$ with essential spectrum $(-\infty,0]$. Its largest eigenvalue $\lambda_0>0$ (the principal eigenvalue) is simple, and the corresponding eigenfunction $\psi$ (the ground state) can be chosen continuous and strictly positive (see, e.g., \cite{RS4}). We normalize it by $\|\psi\|_{L^2(\R^d)} = 1$. Denoting the positions of the particles alive at time $t$ by $X^{(i)}_t$, $1\le i\le \mathcal{N}_t$, the process
\begin{equation} \label{mart}
W_t = e^{-\lambda_0 t}\sum_{i=1}^{\mathcal{N}_t} \psi(X^{(i)}_t)
\end{equation}
is a non-negative martingale, and its limit $W_\infty$ governs the long-time behavior of the whole population. (The event $\{W_\infty = 0\}$ may have positive probability: this is the case if $\beta\not\equiv 0$, since the initial particle may be killed before it branches, and if $d\geq 3$, since it may escape to infinity without ever branching.)

Two natural notions of the population front are the following. The first is the level set $\partial B^x_t = \{y: \rho_1(t,x,y) = 1\}$ of the density. The second is the boundary $\partial A^x_t$ of the region occupied by particles, where $y \in A^x_t$ if the probability of finding at least one particle at time $t$ in a unit neighborhood of $y$ exceeds a fixed value $c$, chosen smaller than $\P(W_\infty>0)$.

Let us compare the structure of the population in branching Brownian motions (BBMs) in homogeneous media, a widely studied topic, with the current case of compactly supported branching. Assume, for simplicity, that $\beta = 0$. In both cases $\partial A^x_t$ and $\partial B^x_t$ grow linearly in $t$, up to lower-order terms, at the same rate.

Inside the front, the number of particles in a domain of fixed size grows exponentially. More precisely, let $U_y$ denote the translation of $U$ by $y\in \R^d$, and let $y = c t\, \bu$ with $|\bu|=1$ and $c$ smaller than the asymptotic speed of the front. Then $n^x_t(U_{ct\bu})$ grows exponentially with an exponent that depends on $c$, and $n^x_t(U_{ct\bu}) / \E n^x_t(U_{ct\bu})$ converges almost surely; in the compactly supported case the limit does not depend on $c$ or $\bu$ (see \cite{Uchiyama} for homogeneous media and \cite{K2} for compactly supported branching). In the compactly supported case, all the moments converge as well \cite{K2}. In the case of constant branching, on the contrary, intermittency is observed: there is $k = k (c)$ such that the $k$-th moment of $n^x_t(U_{ct\bu}) / \E n^x_t(U_{ct\bu})$ tends to infinity as $t \rightarrow \infty$. Such intermittency results have been extended to more general supercritical branching diffusions in periodic media \cite{HKN}. For a fixed set $U$, the asymptotics of $n^x_t(U)$ was described earlier (see, e.g., \cite{W, AH, EHK}), and the convergence of all the moments of the normalized version of $n^x_t(U)$ can be found in \cite{K}; at bounded distances from the origin there is no intermittency in either the constant or the compactly supported case.

We now turn to the subject of this paper, namely the population near the front. In homogeneous media, the evolution of $\partial A^x_t$ is described by the FKPP equation. For $d = 1$ and $\alpha \equiv 1$, Bramson \cite{Bram} showed that $\partial A^x_t$ lags behind $\partial B^x_t$ by a distance that is logarithmic in $t$; this has since been refined and extended to periodic media and to higher dimensions (see, e.g., \cite{N1, N2} and references therein). For $d=1$, let $R_t$ be the position of the rightmost particle at time $t$. Lalley and Sellke \cite{LS87} identified a random variable $Z$ (the limit of the derivative martingale) such that
\begin{equation*}\lim_{t \to \infty} \P\Big(R_t - \sqrt{2}\, t +\frac{3}{2\sqrt{2}}\ln t \leq z \Big)  = \E\exp\{-Z e^{-\sqrt{2} z}\},\qquad z\in\R.\end{equation*}
Moreover, the extremal point process converges in distribution \cite{ABBS, ABK}:
\begin{equation}\label{extremal}\lim_{t\to \infty}\sum_{i = 1}^{\mathcal{N}_t}\delta_{X_t^{(i)} - \sqrt{2} t  + \frac{3}{2\sqrt{2}}\ln t} = \mathcal J, \end{equation}
where $\mathcal{J}$ is a randomly shifted \emph{decorated} Poisson point process. These results were recently extended to higher dimensions \cite{KLZ, BLMZ}.

For compactly supported (or, more generally, rapidly decaying) branching, the picture is different. Erickson \cite{Erik} and Lalley and Sellke \cite{LS1} considered $d=1$; in \cite{LS1} it is shown that $R_t - \sqrt{\lambda_0/2}\, t$ converges in distribution, with
\[\lim_{t \to \infty} \P\big(R_t - \sqrt{\lambda_0/2}\, t \leq z \big)  = \E\exp\{-Z e^{-\sqrt{2\lambda_0} z}\}\]
for a random variable $Z$. For $d\geq 1$, let $R_t$ be the maximal Euclidean norm of the particles alive at time $t$. Shiozawa \cite{S1} showed, for Kato class branching measures with compact support, that $R_t/t \to \sqrt{\lambda_0/2}$ on the event that the population grows; Bocharov and Harris \cite{BH14, BH16} studied the case of a single catalytic point in $d=1$. Nishimori and Shiozawa \cite{NS} identified the logarithmic correction in all dimensions and showed that
\[\lim_{t \to \infty} \P\Big(R_t - \sqrt{\lambda_0/2}\, t  - \frac{d-1}{2\sqrt{2\lambda_0}}\ln t \leq z \Big)  = \E\exp\{-Z e^{-\sqrt{2\lambda_0} z}\},\]
with $Z$ proportional to $W_\infty$.

The distribution of the particles near the frontier in this setting was studied by Bocharov \cite{Boc} for a single catalytic point in $d=1$ and by Nishimori \cite{Nis} for Kato class branching measures with compact support and general offspring distributions. In dimensions $d=1,2$, Nishimori showed that the particles near the frontier, recorded in polar coordinates relative to a sphere of radius $R(t) = \sqrt{\lambda_0/2}\, t + \frac{d-1}{2\sqrt{2\lambda_0}}\ln t + O(1)$, form asymptotically a Poisson point process with random intensity proportional to $W_\infty$. In contrast with (\ref{extremal}), there are no decorations. The test sets in \cite{Nis} have the form $[R(t)+r_1, R(t)+r_2]\times \Theta$ in polar coordinates, with a fixed set $\Theta$ of directions, so that their volume grows like $t^{d-1}$.

In the present paper, we study \emph{test sets of bounded size}, $U_{y(t)}$, located near the front in an arbitrary dimension $d\ge 1$. Since such a window has bounded volume, the relevant distance from the origin is
\[
a(t) = \sqrt{\frac{\lambda_0}{2}}\, t  - \frac{d-1}{2\sqrt{2\lambda_0}}\ln t,
\]
where the logarithmic correction now has the opposite sign; the two corrections differ precisely by the amount that compensates for the factor $t^{d-1}$ in the volume. Our main result (\thmref{mainresult}) states that, if $y(t) = (a(t)-b(t))\bu$ with $b(t)\to b$, then $n^x_t(U_{y(t)})$ converges in distribution, together with all its moments, to a mixed Poisson random variable with random intensity $\kappa_{b,\bu}W_\infty$, where $\kappa_{b,\bu}$ is explicit. 
If $b(t)\to\infty$, then $n^x_t(U_{y(t)})/\E n^x_t(U_{y(t)})$ converges to $W_\infty/\psi(x)$ in $L^p$ for all $p<\infty$; in particular, all its moments converge to those of $W_\infty/\psi(x)$.
In $d=1$ and in the absence of killing, the convergence in distribution in the first statement is contained in \cite{Nis}.
For $d=2$, \cite{Nis} deals with test sets of angular width of order one, and our result is its local counterpart. For $d\geq 3$, although the limiting distribution of the maximal displacement is known \cite{NS}, a mixed Poisson limit for the number of particles near the front does not seem to have been previously established at any spatial scale. We also allow killing and obtain the convergence of all moments, uniformly in the initial point.

The mechanism behind the result is transparent. 

Branching only occurs in $\cD$, and a particle near the front at time $t$ is typically a descendant of a particle that left $\cD$ at time $t/2 + O(\sqrt{t})$ and then traveled without branching, at speed approximately $\sqrt{2\lambda_0}$.
 Conditionally on the history of the process up to an intermediate time, the contributions of different particles are asymptotically independent and rare, which produces a Poisson limit; the random intensity reflects the total ``mass'' accumulated in $\cD$, which is $W_\infty$. 

Our approach differs from that of \cite{Boc, Nis}, which relies on many-to-one and many-to-two lemmas and on second moment estimates. We compute all the correlation functions of the process asymptotically, using the sharp global asymptotics for solutions of parabolic equations with compactly supported potential obtained in \cite{KV}. This yields the convergence of all moments, uniformly in the initial point, and identifies the limit through its moments.

The paper is organized as follows. In Section~\ref{mainres} we formulate the main results and discuss them. Section~\ref{eqmom} collects the equations for the correlation functions and the asymptotics of the fundamental solution. The asymptotics of the correlation functions near the front is obtained in Section~\ref{asympt}, and the main results are proved in Section~\ref{proofmain}.

\section{Main results}\label{mainres}

Throughout the paper, $|\cdot|$ denotes the Euclidean norm, $S^{d-1}$ the unit sphere in $\R^d$ (with $S^0 = \{-1,1\}$), $\dot y = y/|y|$ for $y\neq 0$, and $\langle\cdot,\cdot\rangle$ the inner product. Let $\lambda_0$ and $\psi$ be as in the introduction. For large $|y|$, the ground state satisfies (see \propref{precise} below)
\begin{equation}\label{psi0}
\psi(y) = C_\psi(\dot{y}) |y|^{\frac{1-d}{2}}e^{-\sqrt{2\lambda_0}|y|} (1+O(|y|^{-1})), \qquad |y| \to \infty,
\end{equation}
where $C_\psi$ is a positive continuous function on $S^{d-1}$ (an explicit formula is given in \remref{remC}).

Fix a bounded nonempty open set $U\subset\R^d$ and a unit vector $\bu$. For $t\ge 1$ let
\begin{equation} \label{yfr}
y(t) = \Big(\sqrt{\frac{\lambda_0}{2}}\, t  - \frac{d-1}{2\sqrt{2\lambda_0}}\ln t - b(t)\Big)\bu =  (a(t) - b(t)) \bu,
\end{equation}
where either $b(t)\to b\in\R$ (the test set is near the front), or $b(t) \to  \infty$ and $0\le b(t) \leq  a(t)$ (the window is inside the front). 
If $b(t)\to-\infty$, then $\E n^x_t(U_{y(t)})\to 0$ (by \lemref{fp1}(a) and (\ref{psi0})), and this case is not of interest.

Let
\[
\eta^x(t) = \frac{n^x_t(U_{y(t)})}{\E n^x_t(U_{y(t)})}
\]
be the normalized number of particles in the window. Define
\begin{equation}\label{eqgamma}
\gamma_U(\bu) = \int_U e^{-\sqrt{2\lambda_0} \langle \bu, w\rangle} \, dw, \qquad \kappa_{b,\bu} = \gamma_U(\bu)\, C_\psi(\bu)\, \Big(\frac{\lambda_0}{2}\Big)^{\frac{1-d}{4}} e^{\sqrt{2\lambda_0}\, b}.
\end{equation}
When $U$ is the unit ball, $\gamma_U(\bu) = \int_{|w|<1} e^{\sqrt{2\lambda_0}\, w_1}\,dw$ does not depend on $\bu$. We define the functions $G_k$, $k\in\mathbb{N}$, recursively by $G_1 = \psi$ and
\begin{equation}\label{eqgk}
G_{k}(x) =   \sum_{i = 1}^{k-1}{\binom{k}{i}} \int_{\cD} \int_{0}^{\infty} e^{-k\lambda_0 r}\rho_1(r,x,z)  \alpha(z)  G_i(z)  G_{k-i}(z)\, dr \, dz,\qquad k \geq 2.
\end{equation}
The functions $G_k$ are positive, bounded, and continuous (see \propref{martingale} below). Finally, $S(k,j)$ denotes the Stirling number of the second kind, i.e., the number of ways to partition a set of $k$ elements into $j$ nonempty subsets.

The probabilistic meaning of the functions $G_k$ is given by the following proposition. We write $\P_x$ and $\E_x$ for the probability and expectation when the process starts with one particle at $x$.

\begin{proposition}\label{martingale}
The process $W_t$ defined in (\ref{mart}) is a non-negative $\P_x$-martingale that is bounded in $L^k(\P_x)$ for every $k\geq 1$. Its limit $W_\infty$ satisfies
\[
\E_x W_\infty^k = G_k(x), \qquad k\geq 1,~x\in\R^d.
\]
Moreover, there is $B>0$ such that $0< G_k(x)\le k!\, B^k$ for all $k\ge 1$ and $x\in\R^d$. In particular, $\E_x e^{\theta W_\infty}<\infty$ for $\theta<1/B$, and the distribution of $W_\infty$ under $\P_x$ is uniquely determined by its moments.
\end{proposition}

\begin{theorem}\label{mainresult}
Let $\Gamma\subset\R^d$ be compact, and let $y(t)$ be given by (\ref{yfr}).

(i) Suppose that $b(t) \to b\in\R$. For each $x\in\R^d$, let $N^x_{b,\bu}$ be a random variable which, under $\P_x$ and conditionally on $W_\infty$, has the Poisson distribution with mean $\kappa_{b,\bu} W_\infty$. Then, for each $k\geq 1$,
\begin{equation}\label{factmom}
\lim_{t\to\infty} \E\big[n^x_t(U_{y(t)})\big(n^x_t(U_{y(t)})-1\big)\cdots\big(n^x_t(U_{y(t)})-k+1\big)\big] = \kappa_{b,\bu}^k\, G_k(x),
\end{equation}
\begin{equation}\label{mainform1}
\lim_{t\to \infty}  \E\big(\eta^x(t)^k\big) =  f_{b, \bu}^k(x) := \frac{1}{\psi^k(x)} \sum_{j = 1}^{k}S(k,j)\, G_j(x)\, \kappa_{b,\bu}^{\,j-k},
\end{equation}
uniformly in $x\in\Gamma$. Moreover, $n^x_t(U_{y(t)})$ converges in distribution to $N^x_{b,\bu}$, and $\eta^x(t)$ converges in distribution to $N^x_{b,\bu}/(\kappa_{b,\bu}\psi(x))$; the moments of these limits are $\E_x (N^x_{b,\bu})^k = \sum_{j=1}^k S(k,j)\kappa_{b,\bu}^j G_j(x)$ and $f^k_{b,\bu}(x)$, respectively.

(ii) Suppose that $b(t) \to \infty$ and $0\le b(t)\le a(t)$. Then, for each $k\geq 1$,
\begin{equation}\label{mainform2}
\lim_{t\to \infty}  \E\big(\eta^x(t)^k\big) = f^k(x) := \frac{G_k(x)}{\psi^k(x)}
\end{equation}
uniformly in $x\in\Gamma$. Moreover, $\eta^x(t)\to W_\infty/\psi(x)$ in $L^p(\P_x)$ for every $p<\infty$, and hence in probability.
\end{theorem}

\begin{remark}[Interpretation]\label{interpretation}
(a) The identity (\ref{mainform1}) is the moment formula for a mixed Poisson distribution: if, conditionally on $\Lambda$, $N$ is Poisson with mean $\Lambda$, then $\E N^k = \sum_{j=1}^k S(k,j)\E\Lambda^j$.  As $b\to\infty$, we have $\kappa_{b,\bu}\to\infty$ and $N^x_{b,\bu}/(\kappa_{b,\bu}\psi(x))\to W_\infty/\psi(x)$, which matches part (ii).

(b) The limits in (\ref{mainform1}) and (\ref{mainform2}) are not equal to one, even though the particles near the front are asymptotically ``independent'' in the sense that they descend from different particles that left $\cD$ long ago. Conditionally on the early history of the process they are independent; unconditionally, they share the random factor $W_\infty$. 

(c) At the front, $\eta^x(t)$ does not approach $W_\infty/\psi(x)$: the proof shows that $\lim_{t\to\infty}\E_x\big(\eta^x(t) - W_\infty/\psi(x)\big)^2 = 1/(\kappa_{b,\bu}\psi(x))$ in the setting of part (i), which is exactly the expected conditional Poisson variance of $N^x_{b,\bu}/(\kappa_{b,\bu}\psi(x))$ given $W_\infty$.

(d) In contrast with (\ref{extremal}), the limit contains no decorations: on the scale of bounded windows, the particles near the front do not form clusters.
\end{remark}

\begin{remark}[Relation to \cite{Nis}]\label{general}
In the absence of killing, the martingale $M_t$ of \cite{Nis} coincides with $W_t$, so that the random intensities of the Poisson limits in \cite{Nis} and in \thmref{mainresult} are given by the same random variable $W_\infty$. The proofs of Sections~\ref{asympt} and~\ref{proofmain} do not use the fact that the windows have bounded size: for arbitrary measurable sets $A_t\subset\{w:|w|\le\theta^*t\}$, $\theta^*<\sqrt{2\lambda_0}$, the same arguments show that the $k$-th factorial moment of $n^x_t(A_t)$ is asymptotic to $\big(e^{\lambda_0t}\int_{A_t}\psi(w)\,dw\big)^kG_k(x)$, and similarly for joint factorial moments over several disjoint sets. Applied to the sets considered in \cite{Nis}, this yields the results of \cite[Theorems 2.1--2.3]{Nis} in every dimension, with convergence of all moments and in the presence of killing, albeit under our more restrictive assumptions on the branching mechanism (binary branching with continuous rates, rather than Kato class branching measures and general offspring distributions).
\end{remark}

\begin{remark}[The constant $C_\psi$]\label{remC}
Since $(\lambda_0 - \frac12\Delta)\psi = v\psi$ and $v\psi$ is compactly supported, $\psi$ is the convolution of $v\psi$ with the Green function of $\lambda_0-\frac12\Delta$. Using the asymptotics of the latter (a modified Bessel function), one obtains (\ref{psi0}) with
\[
C_\psi(\bu) = (2\pi)^{-\frac{d-1}{2}}\,(2\lambda_0)^{\frac{d-3}{4}} \int_{\cD} e^{\sqrt{2\lambda_0}\,\langle \bu, z\rangle}\, v(z)\psi(z)\, dz.
\]
The function $C_\psi$ is positive (see Remark 1 in \cite{KV}). 
\end{remark}

\begin{remark}[The two fronts]\label{fronts}
By \propref{precise} below, for $y = (a(t) - b)\bu$ with $b$ in a bounded set,
\[
\rho_1(t,x,y) = \psi(x)\, C_\psi(\bu)\Big(\frac{\lambda_0}{2}\Big)^{\frac{1-d}{4}} e^{\sqrt{2\lambda_0} b}(1+o(1)), \qquad t\to\infty,
\]
uniformly in $\bu\in S^{d-1}$, $x\in\Gamma$. Moreover, by (\ref{ga}) and (\ref{psi0}), $\rho_1(t,x,y)\ge c\,e^{\sqrt{2\lambda_0}B}$ if $|y|\le a(t)-B$, and, by (\ref{psi0}) and the pointwise bound $\rho_1(t,x,y)\le c\,e^{\lambda_0 t}\psi(y)$ obtained in the proof of \lemref{fp1}(a), $\rho_1(t,x,y)\le C e^{-\sqrt{2\lambda_0}B}$ if $|y|\ge a(t)+B$, for large $t$. Hence, in every direction, the front $\partial B^x_t$ lies within a bounded distance of the sphere of radius $a(t)$. As for $\partial A^x_t$, \thmref{mainresult}(i), applied with $U$ the unit ball, gives
\[
\lim_{t\to\infty}\P\big(n^x_t(U_{y(t)}) > 0\big) = 1-\E_x \exp\{-\kappa_{b,\bu}W_\infty\}
\]
for $b(t)\to b$. As $b$ increases from $-\infty$ to $\infty$, the right-hand side increases continuously from $0$ to $\P_x(W_\infty>0)$. Thus, along each ray, the points at distance $a(t)-b$ from the origin belong to $A^x_t$ for all sufficiently large $t$ if $b$ is large enough, and do not belong to $A^x_t$ for all sufficiently large $t$ if $-b$ is large enough. We do not address here the uniformity of this statement with respect to the direction and the depth.
\end{remark}

\begin{remark}[Comparison with \cite{K2}]\label{compareresults}
In \cite{K2}, for $\beta=0$ and $|y(t)| = ct$ with $c < \sqrt{\lambda_0/2}$, the quantity $\tilde{\eta}^x(t) = n^x_t(U_{y(t)})/g(t)$ was considered, where
\[g(t) = e^{\lambda_0 t}\int_{U_{y(t)}} \psi(w)\,dw\, \Big/ \int_{\R^d}\psi(w)\,dw.\]
It was shown that $\tilde{\eta}^x(t)$ converges almost surely and in $L^2$ to a random variable $\tilde{\xi}^x$ with $\E(\tilde{\xi}^x)^k = (\int_{\R^d}\psi(w)\, dw)^k \tilde{f}^k(x)$, where the functions $\tilde f^k$ satisfy exactly the recursion (\ref{eqgk}); i.e., $\tilde f^k = G_k$. By (\ref{comparetogtu}) below, $g(t) = \E n^x_t(U_{y(t)})(1+o(1))/(\psi(x)\int_{\R^d}\psi(w)\,dw)$, so that $\tilde\eta^x(t)$ and $\eta^x(t)$ differ by the deterministic factor $\psi(x)\int_{\R^d}\psi(w)\,dw$ (up to $1+o(1)$). By \propref{martingale}, $\tilde\xi^x$ has the same distribution as $(\int_{\R^d}\psi(w)\,dw)\, W_\infty$.
\end{remark}

\section{Correlation functions and the fundamental solution} \label{eqmom}

\subsection{Correlation functions}
For $t > 0$ and distinct points $y_1, \ldots,y_n \in \R^d$, the $n$-point correlation function $\rho_n(t,x,y_1,\ldots, y_n)$ is defined as the limit, as $\delta\downarrow 0$, of the probability of finding $n$ distinct particles at time $t$ in the balls $B_\delta(y_1) ,\ldots, B_\delta(y_n)$, respectively, divided by the $n$-th power of the volume of  $B_\delta(0)$. (We extend this function by continuity to allow for non-distinct $y_1,\ldots,y_n$.)  For $n=1$, this is the density $\rho_1$, which satisfies (\ref{dens}); as a function of $x$, $\rho_1(t,\cdot,y)$ is the kernel of the semigroup $e^{t\cL}$. For $n \geq 2$ (see \cite[Section 2]{K2}),
\begin{equation}\label{rhon}
\partial_t\rho_n(t,x,y_1,\ldots,y_n) = \cL_x \rho_n(t,x,y_1,\ldots,y_n) + \alpha(x) H_n(t,x,y_1,\ldots,y_n),
\end{equation}
with $\rho_n(0,x,y_1,\ldots,y_n) \equiv 0$,
where
\[
H_n(t,x,y_1,\ldots,y_n) = \sum\limits_{V \subset Y,\, V \neq \emptyset, Y} \rho_{|V|} (t, x, V)\,\rho_{n-|V|}(t, x, Y\setminus V),
\]
$Y = (y_1 , \ldots, y_n )$, the sum is over the proper nonempty subsequences $V$ of $Y$, and $|V|$ is the number of elements in $V$.

For $y\in\R^d$ and $k\geq 1$, define
\[
m_k^y(t,x) = \int_{U_{y}}\cdots\int_{U_{y}} \rho_k(t,x,y_1,\ldots,y_k)\, dy_1\ldots dy_k .
\]
By the definition of the correlation functions, $m_k^y(t,x)$ is the $k$-th factorial moment of $n^x_t(U_y)$:
\[
m_k^y(t,x) = \E\big[n^x_t(U_{y})\big(n^x_t(U_{y})-1\big)\cdots\big(n^x_t(U_{y})-k+1\big)\big].
\]
Consequently,
\begin{equation}\label{kthmoment}
\E\big(n^x_t(U_{y})^k\big)  = \sum_{j = 1}^{k}S(k,j)\, m_j^{y}(t,x).
\end{equation}
Integrating (\ref{rhon}) over $U_y\times\cdots\times U_y$, we obtain
\[
\partial_t m_1^y(t,x) = \cL_xm_1^y(t,x), \qquad m_1^y(0,x) = \chi_{U_{y}}(x),
\]
where $\chi_{U_y}$ is the indicator function of $U_y$, and, for $k \geq 2$,
\begin{equation}\label{mk}
\partial_t m_k^y(t,x) = \cL_xm_k^y(t,x) + \alpha(x) \sum_{i = 1}^{k-1}\binom{k}{i} m_i^y(t,x)m_{k-i}^y(t,x), \qquad m_k^y(0,x) \equiv 0.
\end{equation}
By the Duhamel principle,
\begin{equation}\label{duhamel}
m_k^y(t,x) =  \int_0^t \int_{\cD} \rho_1(t-s,x,z)\,\alpha(z) \sum_{i = 1}^{k-1}\binom{k}{i}m_i^y(s,z)\,m_{k-i}^y(s,z)\, dz\, ds,\qquad k\geq 2.
\end{equation}
We will also need a mixed second-order quantity. Let
\[
m^{y,\psi}(t,x) = \int_{\R^d}\int_{\R^d}\rho_2(t,x,y_1,y_2)\,\chi_{U_y}(y_1)\,\psi(y_2)\, dy_1\, dy_2 .
\]
Since $H_2(t,x,y_1,y_2) = 2\rho_1(t,x,y_1)\rho_1(t,x,y_2)$, the Duhamel form of (\ref{rhon}) with $n=2$ reads $\rho_2(t,x,y_1,y_2) = 2\int_0^t\int_{\cD}\rho_1(t-s,x,z)\alpha(z)\rho_1(s,z,y_1)\rho_1(s,z,y_2)\,dz\,ds$. Integrating it against $\chi_{U_y}(y_1)\psi(y_2)$ (all the integrands are non-negative, so Tonelli's theorem applies) and using \lemref{l33}(a) below, we obtain
\begin{equation}\label{mixedduhamel}
m^{y,\psi}(t,x) = 2\int_0^t\int_{\cD}\rho_1(t-s,x,z)\,\alpha(z)\, m_1^y(s,z)\, e^{\lambda_0 s}\psi(z)\, dz\, ds .
\end{equation}
Moreover,
\begin{equation}\label{mixedmoment}
\E_x\Big[n^x_t(U_y)\sum_{i=1}^{\mathcal N_t}\psi(X^{(i)}_t)\Big] = \int_{U_y}\rho_1(t,x,w)\psi(w)\, dw + m^{y,\psi}(t,x).
\end{equation}

\subsection{Asymptotics of the fundamental solution}
We recall the global asymptotics of $\rho_1$ from \cite{KV}. For $(t,x,y)\in (0,\infty)\times\R^d\times\R^d$, let $\theta = \theta(t,y-x) = |y-x|/t$ and, for $y\neq x$, let $\mathbf{e} = (y-x)/|y-x|$. For $\eps>0$, let
\[\Cint = \{(t,x,y) : \theta\leq \sqrt{2\lambda_0} - \varepsilon\}, \qquad \Cext = \{(t,x,y) : \theta\geq \sqrt{2\lambda_0} +\varepsilon\}.\]
Let $\mathrm{erf}(s) = \frac{2}{\sqrt{\pi}}\int_0^s e^{-r^2}\,dr$, $s\in \R$, and $p_0(t,y) = (2\pi t)^{-\frac{d}{2}}\exp( -|y|^2/2t)$.

\begin{proposition}[\cite{KV}]\label{precise}
Let $\Gamma\subset\R^d$ be compact. The following asymptotic relations hold uniformly in $x\in\Gamma$.
\begin{enumerate}
\item For each $\varepsilon > 0$, there is $\delta>0$ such that
\begin{equation}\label{ga}
\rho_1(t,x,y) = e^{\lambda_0 t} \psi(x)\psi(y) (1+ O(e^{-\delta t}))\quad \text{as } t\to \infty,~ (t,x,y) \in \Cint.
\end{equation}
\item For each $\varepsilon > 0$,
\begin{equation}\label{gaout}
\rho_1(t,x,y) =  p_0(t,y-x) \big(a(\theta,\mathbf{e},x)+ O(t\,|y-x|^{-2})\big)\quad\text{as } |y|\to \infty,~ (t,x,y) \in \Cext,
\end{equation}
where $a$ is a positive continuous function given explicitly in \cite{KV}.
\item For each $\varepsilon_0>0$, as $t\to\infty$ with $\theta\leq 1/\varepsilon_0$,
\begin{equation}\label{intermediate}
\begin{aligned}
  \rho_1(t,x,y) = {}& e^{\lambda_0 t} \psi(x)\psi(y)\Big(1+\mathrm{erf}\Big(\sqrt{t}\,\frac{\sqrt{2\lambda_0}-\theta}{\sqrt{2}}\Big)\Big)\big(a_1(\theta, \dot{y},x)+O(|y|^{-1/2})\big) \\
  &+p_0(t,y-x) \big(a_2(\theta,\mathbf{e},x)+O(|y-x|^{-1})\big),
\end{aligned}
\end{equation}
where $a_1, a_2$ are bounded continuous functions and $a_1(\theta ,\dot{y},x) = 1/2$ when $\theta  \leq  \sqrt{2\lambda_0}$.
\item The ground state satisfies (\ref{psi0}), where $C_\psi$ is a positive continuous function on $S^{d-1}$.
\end{enumerate}
\end{proposition}

\begin{remark}
Part 1 gives a simple expression for $\rho_1$ inside the wide cone $\{\theta\leq \sqrt{2\lambda_0} - \varepsilon\}$, while the region containing the particles, including the front, lies in the narrower cone $\{|y|/t \leq \sqrt{\lambda_0/2}  + \varepsilon \}$. Parts 2 and 3 are only used to obtain the crude upper bound in \lemref{fp1}(a).
\end{remark}

\begin{lemma}\label{l33}
(a) For all $t\ge 0$ and $x\in\R^d$, $\int_{\R^d}\rho_1(t,x,w)\psi(w)\,dw = e^{\lambda_0 t}\psi(x)$. In particular, $\E_x \sum_{i}\psi(X^{(i)}_t) = e^{\lambda_0 t}\psi(x)$.

(b) There is $a_0>0$ such that
\[
\int_{\cD}\rho_1(r,x,z)\,\alpha(z)\, dz \leq a_0e^{\lambda_0 r}, \qquad r\geq 0,~x\in\R^d.
\]
\end{lemma}
\proof
(a) Since $\rho_1(t,\cdot,\cdot)$ is the kernel of $e^{t\cL}$ and $\cL\psi = \lambda_0\psi$, we have
\[
\int_{\R^d}\rho_1(t,x,w)\psi(w)\,dw = (e^{t\cL}\psi)(x) = e^{\lambda_0 t}\psi(x).
\]
The second statement follows since $\E_x \sum_{i}f(X^{(i)}_t) = \int\rho_1(t,x,w)f(w)\,dw$ for bounded measurable $f$.

(b) Since $\psi$ is positive and continuous, $\min_{\cD}\psi>0$, and by (a),
\[
\int_{\cD}\rho_1(r,x,z)\alpha(z)\, dz \le \frac{\|\alpha\|_{\infty}}{\min_{\cD}\psi}\int_{\R^d}\rho_1(r,x,z)\psi(z)\,dz = \frac{\|\alpha\|_{\infty}}{\min_{\cD}\psi}\,e^{\lambda_0 r}\psi(x),
\]
and the claim holds with $a_0 = \|\alpha\|_{\infty}\|\psi\|_\infty/\min_{\cD}\psi$.
\qed

\section{Asymptotics of the correlation functions} \label{asympt}

Throughout this section, $\Gamma$ is a compact set that, without loss of generality, contains $\cD$, and all the asymptotic relations are uniform in $x\in\Gamma$. Let $r_U = \sup_{w\in U}|w|$.

\begin{lemma}\label{windowint}
(a) There is $c_U>0$ such that $\psi(y)\leq c_U\int_{U_y}\psi(w)\,dw$ for all $y\in\R^d$.

(b) As $|y|\to\infty$, uniformly in $\dot y\in S^{d-1}$,
\begin{equation}\label{wint}
\int_{U_y}\psi(w)\, dw = \gamma_U(\dot y)\, C_\psi(\dot y)\, |y|^{\frac{1-d}{2}}e^{-\sqrt{2\lambda_0}|y|}\,(1+o(1)).
\end{equation}

(c) Let $y(t)$ be given by (\ref{yfr}). If $b(t)\to b\in\R$, then $e^{\lambda_0 t}\int_{U_{y(t)}}\psi(w)\, dw\to \kappa_{b,\bu}$. If $b(t)\to\infty$ and $0\le b(t)\le a(t)$, then $e^{\lambda_0 t}\int_{U_{y(t)}}\psi(w)\, dw\to \infty$. If $b(t)\to -\infty$, then $e^{\lambda_0 t}\int_{U_{y(t)}}\psi(w)\, dw\to 0$.
\end{lemma}
\proof
(a) This immediately follows from (\ref{psi0})  and the positivity of $\psi$.

(b) Write $w = y+\tilde w$ with $\tilde w\in U$. As $|y|\to\infty$, uniformly in $\tilde w\in U$ and $\dot y$,
\[
|y+\tilde w| = |y| + \langle\dot y,\tilde w\rangle + O(|y|^{-1}),\qquad \frac{C_\psi(\dot w)}{C_\psi(\dot y)}\to 1,\qquad \frac{|w|^{\frac{1-d}{2}}}{|y|^{\frac{1-d}{2}}}\to 1,
\]
where we used the uniform continuity of $C_\psi$ and the fact that $C_\psi$ is bounded away from zero. Substituting this into (\ref{psi0}), we obtain
\[
\int_{U_y}\psi(w)\, dw = C_\psi(\dot y)|y|^{\frac{1-d}{2}}e^{-\sqrt{2\lambda_0}|y|}\int_U e^{-\sqrt{2\lambda_0}\langle\dot y,\tilde w\rangle}\,d\tilde w\,(1+o(1)),
\]
which is (\ref{wint}).

(c) Note that $\lambda_0 t - \sqrt{2\lambda_0}(a(t)-b(t)) = \frac{d-1}{2}\ln t + \sqrt{2\lambda_0}\, b(t)$. If $b(t)\to b$, then $|y(t)| = a(t)-b(t)$ satisfies $|y(t)|/t\to \sqrt{\lambda_0/2}$, and (\ref{wint}) gives
\[
e^{\lambda_0 t}\int_{U_{y(t)}}\psi(w)\, dw = \gamma_U(\bu)C_\psi(\bu)\Big(\frac{|y(t)|}{t}\Big)^{\frac{1-d}{2}}e^{\sqrt{2\lambda_0}b(t)}(1+o(1))\to \kappa_{b,\bu}.
\]
If $b(t)\to-\infty$, then $|y(t)|\ge a(t)\to\infty$ and $|y(t)|/t\ge \sqrt{\lambda_0/2}\,(1+o(1))$, so that $(|y(t)|/t)^{\frac{1-d}{2}}$ is bounded and the same formula shows that the expression tends to zero. Let $b(t)\to\infty$, $0\le b(t)\le a(t)$, and let $t_n\to\infty$. If $|y(t_n)|$ is bounded along a subsequence, then along it $\int_{U_{y(t_n)}}\psi$ is bounded from below by a positive constant, and the claim is clear. If $|y(t_n)|\to\infty$, then, since $|y(t)|\le a(t)\le \sqrt{\lambda_0/2}\,t$ and $\frac{1-d}{2}\le 0$, we have $(|y(t)|/t)^{\frac{1-d}{2}}\geq (\lambda_0/2)^{\frac{1-d}{4}}$, and (\ref{wint}) gives
\[
e^{\lambda_0 t_n}\int_{U_{y(t_n)}}\psi(w)\, dw\geq \gamma_U(\bu)C_\psi(\bu)\Big(\frac{\lambda_0}{2}\Big)^{\frac{1-d}{4}}e^{\sqrt{2\lambda_0}b(t_n)}(1+o(1))\to\infty.
\]
\qed

\begin{lemma} \label{fp1}
(a) There is $K_1>0$ such that $m_1^y(t,x)\leq K_1 e^{\lambda_0 t} \psi(y)$ for all $t\geq 0$, $y\in\R^d$, and $x\in\Gamma$.

(b) Let $\theta^*\in(0,\sqrt{2\lambda_0})$. There is $\delta>0$ such that, as $t\to\infty$,
\begin{equation}\label{comparetogtu}
m_1^{y}(t,x) =  e^{\lambda_0 t} \psi(x) \int_{U_{y}} \psi(w) \,dw\,  (1+ O(e^{-\delta t}))
\end{equation}
uniformly in $x\in\Gamma$ and $|y|\le\theta^* t$.
\end{lemma}
\proof (a) For $(t,x,w)\in\Cint$, the bound $\rho_1(t,x,w)\le c\, e^{\lambda_0 t}\psi(w)$ follows from (\ref{ga}). In the remaining region, the bound $\rho_1(t,x,w)\le c\, e^{\lambda_0 t}\psi(w)$ for large $t$ follows from (\ref{gaout}), (\ref{intermediate}), (\ref{psi0}), and the elementary inequality $|w-x|^2/(2t)\ge \sqrt{2\lambda_0}|w-x| - \lambda_0 t$, since the Gaussian factor $p_0$ then decays at least as fast as $e^{\lambda_0 t-\sqrt{2\lambda_0}|w|}$, up to factors that are bounded or are dominated by $|w|^{\frac{1-d}{2}}$. For $t$ in a bounded interval, the bound follows from standard Gaussian upper bounds for $\rho_1$. Integrating over $w\in U_y$ and using $\int_{U_y}\psi\le |U|\max_{|w-y|\le r_U}\psi\le c\,\psi(y)$ (by the argument in the proof of \lemref{windowint}(a)), we obtain~(a).

(b) If $|y|\le\theta^*t$, $w\in U_y$, and $x\in\Gamma$, then $|w-x|/t\le \theta^* + (r_U + \max_{\Gamma}|x|)/t\le\sqrt{2\lambda_0}-\eps$ with $\eps = (\sqrt{2\lambda_0}-\theta^*)/2$, for all sufficiently large $t$. It remains to integrate (\ref{ga}) over $w\in U_y$.
\qed

\begin{lemma}\label{fp3} For each $k \geq 1$, the following statements hold.

(a) There is $K_k >0$ such that $m_k^y(t,x)\leq K_k e^{k\lambda_0 t} \psi^k(y)$ for all $t\ge 0$, $y\in\R^d$, and $x\in\Gamma$.

(b) For each $\theta^*\in(0,\sqrt{2\lambda_0})$, there is $\delta_k > 0$ such that, as $t \rightarrow \infty$,
\begin{equation}\label{twt}
m_k^{y}(t,x) =
 e^{k\lambda_0 t}\Big(\int_{U_{y}} \psi(w) \,dw\Big)^k G_k(x)  (1+ O(e^{-\delta_k t }))
\end{equation}
uniformly in $x\in\Gamma$ and $|y|\le\theta^*t$.
\end{lemma}
\proof
We use induction on $k$. For $k = 1$, both statements are contained in \lemref{fp1}. Let $k\geq 2$, and assume that (a) and (b) hold for all $i < k$ (in (b), for all $\theta^*\in(0,\sqrt{2\lambda_0})$).

(a) By (\ref{duhamel}), the induction hypothesis (recall that $\cD\subset\Gamma$), and \lemref{l33}(b),
\[
m_k^y(t,x) \leq  \psi^{k}(y) \sum_{i = 1}^{k-1}\binom{k}{i} K_iK_{k-i}\int_0^t e^{k\lambda_0 s}\int_{\cD} \rho_1(t-s,x,z)\alpha(z)\, dz \, ds \leq K_k e^{k\lambda_0 t}\psi^{k}(y),
\]
with $K_k = \frac{a_0}{(k-1)\lambda_0}\sum_{i = 1}^{k-1}\binom{k}{i} K_iK_{k-i}$.

(b) Fix $\theta^*\in(0,\sqrt{2\lambda_0})$, choose $\theta^{**}\in(\theta^*,\sqrt{2\lambda_0})$, and choose $\eta\in(0,1/2)$ such that $\theta^*\le (1-\eta)\theta^{**}$. We split the time integral in (\ref{duhamel}) into the integrals $J_1$ over $[0,(1-\eta)t]$ and $J_2$ over $[(1-\eta)t, t]$.

For $s\in[0,(1-\eta)t]$, we use the bounds from (a) for $i<k$ and \lemref{l33}(b):
\[
\begin{split}
J_1 &\le a_0\, \psi^{k}(y)\sum_{i = 1}^{k-1}\binom{k}{i} K_i K_{k-i}\int_0^{(1-\eta)t}e^{\lambda_0(t-s)} e^{k\lambda_0 s}\,ds \\
&\le c'\psi^{k}(y)\,e^{k\lambda_0 t}e^{-(k-1)\lambda_0\eta t}.
\end{split}
\]
By \lemref{windowint}(a) and since $G_k$ is bounded from below on $\Gamma$ by a positive constant (\propref{martingale}),
\[
J_1 = e^{k\lambda_0 t}\Big(\int_{U_{y}} \psi(w) \,dw\Big)^k G_k(x)\cdot O(e^{-(k-1)\lambda_0\eta t}).
\]

For $s\in[(1-\eta)t, t]$ and $|y|\le\theta^*t$, we have $|y|\le (1-\eta)\theta^{**}t\le\theta^{**}s$. By the induction hypothesis (b), applied with $\theta^{**}$ in place of $\theta^*$ and with the starting point $z\in\cD\subset\Gamma$, we have $m_i^y(s,z) = e^{i\lambda_0 s}(\int_{U_y}\psi)^i G_i(z)(1+O(e^{-\delta' t}))$ for $i<k$ and some $\delta'>0$. Since the integrand is non-negative, after the change of variables $r = t-s$ we obtain
\[
\begin{split}
J_2 = {}& e^{k\lambda_0 t}\Big(\int_{U_{y}} \psi(w) \,dw\Big)^k (1+ O(e^{-\delta' t}))\\
&\times\int_0^{\eta t}\int_{\cD} e^{-k\lambda_0 r}\rho_1(r,x,z)\alpha(z)\sum_{i = 1}^{k-1}\binom{k}{i}  G_i(z)  G_{k-i}(z)\, dz \, dr.
\end{split}
\]
By \lemref{l33}(b) and the boundedness of the functions $G_i$,
\[
\begin{split}
\int_{\eta t}^\infty\int_{\cD} e^{-k\lambda_0 r}\rho_1(r,x,z)\alpha(z)G_i(z)G_{k-i}(z)\, dz \, dr &\leq a_0\|G_i\|_\infty\|G_{k-i}\|_\infty\int_{\eta t}^\infty e^{-(k-1)\lambda_0 r}dr \\
&= O(e^{-(k-1)\lambda_0\eta t}).
\end{split}
\]
Therefore, the last integral in the formula for $J_2$ equals $G_k(x)(1+O(e^{-(k-1)\lambda_0\eta t}))$. Combining the asymptotics of $J_1$ and $J_2$, we obtain (\ref{twt}) with $\delta_k = \min(\delta', (k-1)\lambda_0\eta)$.
\qed

\begin{lemma}\label{mixed}
For each $\theta^*\in(0,\sqrt{2\lambda_0})$, there is $\delta>0$ such that, as $t\to\infty$,
\[
m^{y,\psi}(t,x) = e^{2\lambda_0 t}\Big(\int_{U_{y}} \psi(w) \,dw\Big) G_2(x)  (1+ O(e^{-\delta t }))
\]
uniformly in $x\in\Gamma$ and $|y|\le\theta^*t$.
\end{lemma}
\proof
The formula (\ref{mixedduhamel}) coincides with (\ref{duhamel}) for $k=2$, except that one of the two factors $m_1^y(s,z)$ is replaced by $e^{\lambda_0 s}\psi(z) = e^{\lambda_0 s}G_1(z)$. The proof of \lemref{fp3}(b) for $k=2$ applies without changes: the contribution of $s\in[0,(1-\eta)t]$ is bounded by $c\,\psi(y)e^{(2-\eta)\lambda_0 t}$, which is negligible by \lemref{windowint}(a), and the contribution of $s\in[(1-\eta)t,t]$ equals
\[
e^{2\lambda_0 t}\Big(\int_{U_{y}} \psi(w) \,dw\Big)(1+O(e^{-\delta' t}))\, 2\int_0^{\eta t}\int_{\cD}e^{-2\lambda_0 r}\rho_1(r,x,z)\alpha(z)\psi^2(z)\, dz\, dr,
\]
where the last integral equals $G_2(x)(1+O(e^{-\lambda_0\eta t}))$.
\qed

\section{Proofs of the main results}\label{proofmain}

\subsection{Proof of \propref{martingale}}
Let $u_k(t,x) = \E_x\big[\big(\sum_{i}\psi(X^{(i)}_t)\big)^k\big]$, $k\geq 0$, so that $u_0\equiv 1$. These quantities are finite, since $\psi$ is bounded and all the moments of the total number of particles at time $t$ are finite. Conditioning on the first event that happens to the initial particle, which is either branching at rate $\alpha$ (after which the population consists of two independent copies of the process started at the branching point) or killing at rate $\beta$ (after which the sum is zero), we obtain
\[
\partial_t u_k = \frac{1}{2}\Delta u_k + \alpha\Big(\sum_{i=0}^{k}\binom{k}{i}u_iu_{k-i} - u_k\Big) - \beta u_k,\qquad u_k(0,x) = \psi^k(x).
\]
Since the terms with $i=0$ and $i=k$ in the sum equal $u_k$, this can be rewritten as
\[
\partial_t u_k = \cL u_k + \alpha\sum_{i=1}^{k-1}\binom{k}{i}u_iu_{k-i},\qquad u_k(0,x) = \psi^k(x),
\]
and therefore
\begin{equation}\label{uk}
\begin{split}
u_k(t,x) = {}&\int_{\R^d}\rho_1(t,x,w)\psi^k(w)\,dw \\
&+ \sum_{i=1}^{k-1}\binom{k}{i}\int_0^t\int_{\cD}\rho_1(t-s,x,z)\alpha(z)u_i(s,z)u_{k-i}(s,z)\,dz\,ds .
\end{split}
\end{equation}

\emph{Martingale property and moment bounds.} By \lemref{l33}(a) and the branching property, $\E_x[W_{t+s}|\cF_t] = e^{-\lambda_0(t+s)}\sum_i e^{\lambda_0 s}\psi(X^{(i)}_t) = W_t$. We claim that $u_k(t,x)\le L_k e^{k\lambda_0 t}$ for all $t\ge 0$ and $x\in\R^d$, with $L_1 = \|\psi\|_\infty$ and $L_k = \|\psi\|^k_\infty + \frac{a_0}{(k-1)\lambda_0}\sum_{i=1}^{k-1}\binom{k}{i}L_iL_{k-i}$. Indeed, $u_1(t,x) = e^{\lambda_0 t}\psi(x)$ by \lemref{l33}(a). For $k\ge 2$, the first term in (\ref{uk}) is at most $\|\psi\|_\infty^{k-1}e^{\lambda_0 t}\psi(x)\le\|\psi\|_\infty^k e^{k\lambda_0 t}$, and, by the induction hypothesis and \lemref{l33}(b), the second term is at most
\[
\sum_{i=1}^{k-1}\binom{k}{i}L_iL_{k-i}\int_0^te^{k\lambda_0 s}a_0e^{\lambda_0(t-s)}\,ds\le \frac{a_0}{(k-1)\lambda_0}\sum_{i=1}^{k-1}\binom{k}{i}L_iL_{k-i}\,e^{k\lambda_0 t}.
\]
Thus $\sup_t\E_xW_t^k\le L_k$, so the martingale $W_t$ is bounded in every $L^k(\P_x)$ and converges almost surely and in every $L^k(\P_x)$ to $W_\infty$. In particular, $\E_xW_\infty^k = \lim_{t\to\infty}g_k(t,x)$, where $g_k(t,x) = e^{-k\lambda_0 t}u_k(t,x)$.

\emph{Identification of the moments.} Multiplying (\ref{uk}) by $e^{-k\lambda_0 t}$ and substituting $r = t-s$, we obtain
\[
\begin{split}
g_k(t,x) = {}& e^{-k\lambda_0 t}\int_{\R^d}\rho_1(t,x,w)\psi^k(w)\,dw \\
&+ \sum_{i=1}^{k-1}\binom{k}{i}\int_0^t\int_{\cD}e^{-k\lambda_0 r}\rho_1(r,x,z)\alpha(z)g_i(t-r,z)g_{k-i}(t-r,z)\,dz\,dr .
\end{split}
\]
We show by induction that $g_k(t,x)\to G_k(x)$ for every $x\in\R^d$. For $k=1$, $g_1(t,x) = \psi(x) = G_1(x)$. Let $k\geq 2$. The first term is at most $\|\psi\|_\infty^ke^{-(k-1)\lambda_0 t}\to 0$. In the second term, $g_i(t-r,z)\to G_i(z)$ as $t\to\infty$ for fixed $r$ and $z$ by the induction hypothesis, and the integrand is dominated by $L_iL_{k-i}e^{-k\lambda_0 r}\rho_1(r,x,z)\alpha(z)$, which is integrable over $(r,z)\in[0,\infty)\times\cD$ by \lemref{l33}(b). By dominated convergence, $g_k(t,x)\to G_k(x)$, which proves that $\E_xW_\infty^k = G_k(x)$.

\emph{Bounds on $G_k$.} Positivity of $G_k$ follows from (\ref{eqgk}) since $\rho_1>0$, $\alpha\geq 0$, $\alpha\not\equiv 0$, and $G_1 = \psi>0$; 
continuity follows from (\ref{eqgk}), since the resolvent $(n\lambda_0-\cL)^{-1}$, $n\ge 2$, maps bounded measurable functions into continuous ones.

 For $n\geq 2$, consider the operator on bounded measurable functions
\[T_n(f)(x) = \int_{\cD} \int_{0}^{\infty} e^{-n\lambda_0 r}\rho_1(r,x,z)  \alpha(z) f(z) \, dr \, dz.\]
By \lemref{l33}(b), $\|T_n\| \leq a_0/((n-1)\lambda_0)$.  Let $A = \max(1, 2a_0/\lambda_0)$. Notice that $\|T_n\| \leq A/n$. 
We claim that
\begin{equation}\label{induchyp}
\|G_n\|_\infty \leq  A^{2n-1}\, n!\, \|\psi\|_\infty^n,\qquad n\geq 1.
\end{equation}
For $n=1$ this holds since $A\geq 1$. If it holds for all indices smaller than $n\geq 2$, then, since $G_n = \sum_{i=1}^{n-1}\binom{n}{i}T_n(G_iG_{n-i})$,
\[
\begin{split}
\|G_{n}\|_{\infty} &\leq \frac{A}{n}\sum_{i = 1}^{n-1}\binom{n}{i} A^{2i-1} i!\,A^{2(n-i)-1} (n-i)!\,\|\psi\|_\infty^n \\
&= A^{2n-1} (n-1)!\,(n-1)\,\|\psi\|_\infty^n\le A^{2n-1}\, n!\,\|\psi\|_\infty^n .
\end{split}
\]
Thus $G_k(x)\le k!B^k$ with $B = A^2\|\psi\|_\infty$. Consequently, $\E_xe^{\theta W_\infty} = \sum_k\theta^kG_k(x)/k!<\infty$ for $0\le\theta<1/B$, and the distribution of $W_\infty$ is determined by its moments (see, e.g., \cite[Section 30]{Bil}).
\qed

\subsection{Proof of \thmref{mainresult}(i)}
Without loss of generality, $\cD\subset\Gamma$. Let $\theta^* = \frac{3}{4}\sqrt{2\lambda_0}$. Since $b(t)$ is bounded and $a(t)\le\sqrt{\lambda_0/2}\,t = \frac12\sqrt{2\lambda_0}\,t$, we have $|y(t)|\le\theta^*t$ for all sufficiently large $t$. By \lemref{fp3}(b) and \lemref{windowint}(c), for each $j\geq 1$,
\[
m_j^{y(t)}(t,x) = \Big(e^{\lambda_0 t}\int_{U_{y(t)}}\psi(w)\,dw\Big)^jG_j(x)(1+o(1))\to\kappa_{b,\bu}^jG_j(x)
\]
uniformly in $x\in\Gamma$, which is (\ref{factmom}). By (\ref{kthmoment}) and (\ref{factmom}),
\[
\lim_{t\to\infty}\E\big(n^x_t(U_{y(t)})^k\big) = \sum_{j=1}^kS(k,j)\kappa_{b,\bu}^jG_j(x)
\]
uniformly in $x\in\Gamma$. Since $\E n^x_t(U_{y(t)}) = m_1^{y(t)}(t,x)\to\kappa_{b,\bu}\psi(x)$ and $\psi$ is bounded from below on $\Gamma$, we obtain (\ref{mainform1}).

If $N$ is Poisson distributed with mean $\Lambda$ conditionally on $\Lambda\ge 0$, then its $j$-th factorial moment equals $\E\Lambda^j$. By \propref{martingale}, the $j$-th factorial moment of $N^x_{b,\bu}$ is therefore $\kappa_{b,\bu}^jG_j(x)$, which, by (\ref{factmom}), is the limit of the $j$-th factorial moment of $n^x_t(U_{y(t)})$. By (\ref{kthmoment}), all the moments of $n^x_t(U_{y(t)})$ converge to those of $N^x_{b,\bu}$.

By \propref{martingale}, for $\theta>0$ small enough that $\kappa_{b,\bu}(e^\theta-1)<1/B$,
\[
\E_x e^{\theta N^x_{b,\bu}} = \E_x\exp\{\kappa_{b,\bu}(e^\theta-1)W_\infty\}<\infty,
\]
so the distribution of $N^x_{b,\bu}$ is determined by its moments. By the method of moments (see \cite{FS} or \cite[Section 30]{Bil}), $n^x_t(U_{y(t)})$ converges in distribution to $N^x_{b,\bu}$. Since $\E n^x_t(U_{y(t)})\to\kappa_{b,\bu}\psi(x)$, it follows that $\eta^x(t)$ converges in distribution to $N^x_{b,\bu}/(\kappa_{b,\bu}\psi(x))$.
\qed

\subsection{Proof of \thmref{mainresult}(ii)}
Again, let $\cD\subset\Gamma$ and $\theta^* = \frac34\sqrt{2\lambda_0}$, so that $|y(t)|\le a(t)\le\theta^*t$. Let $\Phi(t) = e^{\lambda_0 t}\int_{U_{y(t)}}\psi(w)\,dw$. By \lemref{windowint}(c), $\Phi(t)\to\infty$, and by \lemref{fp3}(b), $m_j^{y(t)}(t,x) = \Phi(t)^jG_j(x)(1+o(1))$ uniformly in $x\in\Gamma$. Since the functions $G_j$ are bounded and bounded from below on $\Gamma$ by positive constants, the terms with $j<k$ in (\ref{kthmoment}) are negligible compared with the term $j=k$, and
\[
\E\big(n^x_t(U_{y(t)})^k\big) = \Phi(t)^kG_k(x)(1+o(1)),\qquad \E n^x_t(U_{y(t)}) = \Phi(t)\psi(x)(1+o(1)).
\]
This implies (\ref{mainform2}) uniformly in $x\in\Gamma$.

To prove the $L^2$ convergence, we write
\[
\E_x\Big(\eta^x(t)-\frac{W_\infty}{\psi(x)}\Big)^2 = \E_x\big(\eta^x(t)\big)^2 - \frac{2}{\psi(x)}\E_x\big[\eta^x(t)W_\infty\big] + \frac{\E_xW_\infty^2}{\psi^2(x)}.
\]
The first term tends to $G_2(x)/\psi^2(x)$ by (\ref{mainform2}), and the last term equals $G_2(x)/\psi^2(x)$ by \propref{martingale}. Since $W_t$ is a martingale bounded in $L^2$ and $n^x_t(U_{y(t)})$ is $\cF_t$-measurable, $\E_x[n^x_t(U_{y(t)})W_\infty] = \E_x[n^x_t(U_{y(t)})W_t]$, and, by (\ref{mixedmoment}),
\[
\E_x[n^x_t(U_{y(t)})W_t] = e^{-\lambda_0 t}\int_{U_{y(t)}}\rho_1(t,x,w)\psi(w)\,dw + e^{-\lambda_0 t}m^{y(t),\psi}(t,x).
\]
The first term on the right-hand side is at most $\|\psi\|_\infty e^{-\lambda_0 t}m_1^{y(t)}(t,x)$, and, by \lemref{mixed}, the second term equals $\Phi(t)G_2(x)(1+o(1))$. Dividing by $\E n^x_t(U_{y(t)}) = \Phi(t)\psi(x)(1+o(1))$, we obtain $\E_x[\eta^x(t)W_\infty]\to G_2(x)/\psi(x)$. Therefore,
\[
\lim_{t\to\infty}\E_x\Big(\eta^x(t)-\frac{W_\infty}{\psi(x)}\Big)^2 = \frac{G_2(x)}{\psi^2(x)} - \frac{2G_2(x)}{\psi^2(x)} + \frac{G_2(x)}{\psi^2(x)} = 0 .
\]
Thus $\eta^x(t)\to W_\infty/\psi(x)$ in $L^2(\P_x)$ and in probability. Let $p<\infty$ and $k = 2\lceil p\rceil$. By (\ref{mainform2}), $\sup_{t\ge t_0}\E_x(\eta^x(t))^{k}<\infty$, and $\E_x W_\infty^k<\infty$ by \propref{martingale}. Hence the family $|\eta^x(t) - W_\infty/\psi(x)|^p$, $t\geq t_0$, is bounded in $L^{k/p}$ with $k/p>1$, and therefore uniformly integrable. Consequently, the convergence holds in $L^p(\P_x)$.
\qed

\begin{remark}
The same computation in the setting of part (i) gives $\E_x(\eta^x(t))^2\to G_2(x)/\psi^2(x)+1/(\kappa_{b,\bu}\psi(x))$, while the other two terms have the same limits as above (note that \lemref{mixed} applies since $|y(t)|\le\theta^*t$). Hence $\E_x(\eta^x(t)-W_\infty/\psi(x))^2\to 1/(\kappa_{b,\bu}\psi(x))$, as stated in \remref{interpretation}(c).
\end{remark}

\noindent\textbf{Acknowledgments:} Pratima Hebbar was supported by the pre-tenure research leave granted by Grinnell College.
Leonid Koralov was supported by the NSF grant DMS-2307377 and the Simons Foundation Grant MP-TSM-00002743.

\end{document}